\documentclass[12pt,a4paper]{article}
\usepackage{amsmath, amsfonts, amssymb,amsthm}

\begin{document}
\begin{center}
APPROXIMATION by POLYNOMIALS in a WIEGHTED SPACE of

INFINITELY DIFFERENTIABLE FUNCTIONS.
\end{center}
\begin{center}
P. V. FEDOTOVA, I. Kh. Musin
\end{center}
\newtheorem{theorem}{Theorem}
\newtheorem*{lemma}{Lemma}

\begin{abstract}
{\sc Abstract}. The density of polynomials in a weighted space of
infinitely differentiable functions on ${\mathbb R}^n$ is proved under
minimal conditions on weight functions and on differences between weight
functions. We apply this result for description of strong dual for weighted
spaces of infinitely differentiable functions on ${\mathbb R}$ and weighted
spaces of sequences of infinitely differentiable functions ${\mathbb
R}$ in terms of the Fourier-Laplace transform of functionals.
\end{abstract}

\vspace {1cm}

1. {\it On approximation by polynomials}.
Let $\varphi=\{\varphi_m\}_{m=1}^{\infty}$ be a family of continuous
functions $\varphi_m: {\mathbb R}^n \to {\mathbb R}$ such that:

1). $
\lim \limits_{x \to \infty}
\displaystyle \frac {\varphi_m(x)} {\Vert x \Vert} = +\infty $
for each $m \in {\mathbb N}$;

2). $\lim \limits_{x \to \infty}
(\varphi_m(x) - \varphi_{m+1}(x)) = +\infty , \
 m \in {\mathbb N}$.

For each $m \in {\mathbb N}$ \
let
$$
{\cal E}(\varphi_m) = \{f \in C^m({\mathbb R}^n):
q_m (f) = \sup_{x \in {\mathbb R}^n, \vert \alpha \vert \le m}
\displaystyle \frac
{\vert D^{\alpha} f(x) \vert}
{\exp(\varphi_m(x))} <  \infty \}.
$$
Obviously,
${\cal E}(\varphi_{m+1}) \subset {\cal E}(\varphi_m)$.

Set
$
{\cal E}(\varphi) =
\displaystyle \bigcap \limits_{m=1}^{\infty} {\cal E}(\varphi_m).
$
Endow ${\cal E}(\varphi)$ with the topology of projective limit of the spaces
${\cal E}(\varphi_m)$.
Clearly, the topology of ${\cal E}(\varphi)$ can be defined by the family of
norms
$$
q_{p,m}(f) =
\sup_{x \in {\mathbb R}^n, \vert \alpha \vert \le p}
\displaystyle \frac
{\vert D^{\alpha} f(x) \vert}
{\exp(\varphi_m(x))} \ ,
\
f \in {\cal E}(\varphi), p \in {\mathbb Z_+}, m \in {\mathbb N}.
$$

Note that for each $m \in {\mathbb N}$ the space
${\cal E}(\varphi_{m+1})$ is completely continuously
embedded into the space ${\cal E}(\varphi_m)$ [1].
So ${\cal E}(\varphi)$ is an $(M^*)$-space [2].

Under additional assumptions on weight functions $\varphi_m$ the space
${\cal E}(\varphi)$ was considered from various points of view by L.
Ehrenpreis [3], B.A. Taylor [4], I. Kh. Musin [5], [6] and others.
In particular, the problem of polynomial approximation in ${\cal E}(\varphi)$
was studied in [5], [6] in the case when $\varphi_m(x) = \varphi (x) -
m \ln (1 + \Vert x \Vert)$, $\varphi$ is a real function on ${\mathbb
R}^n$ such that for some $C>0, D \in {\mathbb R}$ \ $\varphi (x) \ge
C{\Vert x \Vert}^{\mu} - D, \ x \in {\mathbb R}^n.$
Now under minimal conditions on weight functions and on differences between
weight functions we study this problem again.

Set for compactness
$\theta_m(x) = \exp (\varphi_m(x)), \
x \in {\mathbb R}^n.
$

For a function $u: [0, \infty) \to {\mathbb R}$ such that
\begin{equation}
\displaystyle \lim \limits_{x \to +\infty}
\displaystyle \frac {u(x)} {x} = + \infty
\end{equation}
let
$
u^*(x) = \sup \limits_{y \ge 0}(x y - u(y)), \
u(e)(x) = u(e^x), \
x \ge 0.
$
Note that $u^*(x) < \infty $ on $[0, \infty)$ and
$
\displaystyle \lim \limits_{x \to +\infty}
\displaystyle \frac {u^*(x)} {x} = + \infty .
$

The following useful lemma holds. Its proof is elementary.
\begin{lemma}
Let a lower semi-continuous function $u: [0, \infty) \to {\mathbb R}$
satisfy (1). Then
$$
(u(e))^*(x) + (u^*(e))^*(x) \le x \ln x - x, \
x > 0.
$$
\end{lemma}

\begin{theorem}
The polynomials are dense in $ {\cal E}(\varphi) $.
\end{theorem}

{\bf Proof}. The beginnig of the proof is the same as the one in the
proof of Lemma 6 in [6].

Let $ f \in {\cal E}(\varphi) $, that is, $f
\in {\cal E}({\mathbb R}^n)$ and for each $m \in {\mathbb N}$ there
exists $c_m >0$ such that
\begin{equation}
\vert D^{\alpha} f(x) \vert
\le c_m \theta_m(x), \ x \in {\mathbb R}^n, \ \vert \alpha \vert \le m.
\end{equation}

  Let us approximate $f$ by polynomials in ${\cal E}(\varphi)$.
There are three steps in the proof.

1. For $r >0$ let
$\Pi_r=\{x \in {\mathbb R}^n: \vert x_j \vert < r, j =1, \ldots , n \}$.
We consider a function $\chi \in C^{\infty}({\mathbb R})$  such that
supp  $ \chi \subseteq [-2,2],
\chi (x) = ~1 $ for $ x \in [-1,1]$, $ 0 \le  \chi (x) \le 1 $
$ \forall x \in {\mathbb R}$.
We set
$
\eta (x_1, x_2, \ldots , x_n) = \chi (x_1) \chi (x_2) \cdots \chi (x_n).
$
Let
$ f_{\nu}(x) = f(x)\eta (\frac x {\nu}), \ {\nu} \in {\mathbb N},
x \in {\mathbb R}^n $.
Obviously, $f_{\nu} \in {\cal E}(\varphi)$.
We claim that
$f_{\nu}\to $ $f$ in ${\cal E}(\varphi)$ as ${\nu} \to \infty$.
For each
$m \in {\mathbb N}$ we have
$$
\sup \limits_{x \in {\mathbb R}^n} \displaystyle \frac
{\vert f_{\nu}(x) - f(x) \vert} {\theta_m(x)} \le
\sup \limits_{x \notin \Pi_{\nu}}
\displaystyle \frac {\vert f(x)
\vert}{{\theta}_m(x)}
\le
\sup \limits_{x \notin \Pi_{\nu}}
\displaystyle \frac {c_{m+1} {\theta}_{m+1}(x)}
{{\theta}_m(x)}.
$$
Hence, as
${\nu} \to \infty$
\begin{equation}
\sup \limits_{x \in {\mathbb R}^n} \displaystyle \frac
{\vert f_{\nu}(x) - f(x) \vert} {\theta_m(x)} \to 0.
\end{equation}
Now \footnote{
For two multi-indices
$\alpha = ({\alpha}_1, \ldots , {\alpha}_n),
\beta = ({\beta}_1, \ldots , {\beta}_n) \in {\mathbb Z_+^n}$
the notation
$\beta \le \alpha $ indicates that
${\beta}_j \le {\alpha}_j, \ j = 1, 2, \ldots , n$.},
$$
\sup \limits_{x \in {\mathbb R}^n,
1 \le \vert \alpha \vert \le m}
\displaystyle \frac
{\vert D^{\alpha} (f_{\nu}(x) - f(x)) \vert} {\theta_m(x)} =
$$
$$
\sup \limits_{x \in {\mathbb R}^n,
1 \le \vert \alpha \vert \le m}
\displaystyle \frac
{\vert \sum \limits_{\beta \le \alpha, \vert \beta \vert < \vert \alpha \vert}
\left(\beta \atop \alpha \right)
(D^{\beta} f)(x){\nu}^{\vert \beta \vert - \vert \alpha \vert}
(D^{\alpha - \beta} \eta)(\frac x {\nu})
+
(D^{\alpha}f)(x)(\eta (\frac x {\nu}) - 1) \vert}
{{\theta}_m(x)}
$$
$$
\le
\sup \limits_{x \in \Pi_{2\nu}
\setminus \Pi_{\nu}, 1 \le \vert \alpha \vert \le m}
\displaystyle \frac
{\sum \limits_{\beta \le \alpha,
\vert \beta \vert < \vert \alpha \vert}
\left(\beta \atop \alpha \right)
{\nu}^{\vert \beta \vert - \vert \alpha \vert}
\vert (D^{\beta} f)(x) \vert
\vert (D^{\alpha - \beta} \eta)(\frac x {\nu})\vert}
{{\theta}_m(x)}
+
$$
$$
+ \sup \limits_{x \notin \Pi_{\nu}, 1 \le \vert \alpha \vert \le m}
\displaystyle \frac
{\vert (D^{\alpha}f)(x) \vert}
{{\theta}_m(x)}
$$
Hence, using inequality (2) we conclude that as $ {\nu} \to \infty$
$$
\sup \limits_{x \in {\mathbb R}^n,
1 \le \vert \alpha \vert \le m}
\displaystyle \frac
{\vert
D^{\alpha}
(f_{\nu}(x) - f(x)) \vert} {\theta_m(x)} \to 0.
$$
It now follows by (3) that
$ q_m(f_{\nu} - f) \to 0 $
as $ \nu \to \infty $.
Since $ m \in {\mathbb N}$ can be arbitrary this means that the sequence
$ (f_{\nu})_{\nu=1}^{\infty} $
converges to $f$ in ${\cal E}(\varphi)$ as $ \nu \to \infty $.

2. Fix   $ \nu \in {\mathbb N}$. Let $h$ be an entire function (not
identical zero) of exponential type at most 1 such that $ h \in
L_1({\mathbb R}),  h(x) \ge 0, x \in  {\mathbb R} $. For example, we
can take $ h(z) = \displaystyle \frac  {\sin^2{ \frac z 2}} {z ^2}, \
z \in~{\mathbb C} $. Let $H(z_1, z_2, \ldots , z_n) = h(z_1) h(z_2)
\cdots h(z_n)$. By the Paley-Wiener theorem ([7], Chapter 6) we can
find a positive constant $C_H >0$ such that for each $\alpha \in
{\mathbb Z_+}^n$
\begin{equation}
\vert
(D^{\alpha} H)(x) \vert \le C_H \ , \ x \in {\mathbb R}^n.
\end{equation}

Let
$ \int_{{\mathbb R}^n} H(x) \ d x = A $.
For $\lambda> 1$ we now set
$$
f_{\nu,\lambda}(x) =  \displaystyle \frac  {{\lambda}^n} A
\int_{{\mathbb R}^n} f_{\nu}(y)
H(\lambda (x-y)) \ d y, \ x \in {\mathbb R}^n.
$$
Obviously,
$f_{\nu,\lambda} \in {\cal E}(\varphi)$.
Let us show that
$f_{\nu,\lambda} \to f_{\nu} $ in ${\cal E}(\varphi)$
as
$ \lambda \to +\infty $.

We consider arbitrary $m \in {\mathbb N}$ and let
$r(\lambda) = {\lambda}^{-\frac {2n}{2n+1}}$.
For arbitrary $ \alpha \in {\mathbb Z}_+, x \in {\mathbb R}^n$
$$
(D^{\alpha}f_{\nu, \lambda})(x) -
(D^{\alpha}f_{\nu})(x)
=
\displaystyle \frac  {{\lambda}^n} A
\int_{{\mathbb R}^n}
((D^{\alpha}f_{\nu})(y) - (D^{\alpha}f_{\nu})(x))
H(\lambda (x-y)) \ dy
$$
$$
= \displaystyle \frac  {{\lambda}^n} A
\int_{\Vert y - x \Vert \le r(\lambda)}
((D^{\alpha}f_{\nu})(y) - (D^{\alpha}f_{\nu})(x))
H(\lambda (x-y)) \ dy
$$
$$
+ \displaystyle \frac  {{\lambda}^n} A
\int_{\Vert y - x \Vert > r(\lambda)}
((D^{\alpha}f_{\nu})(y) - (D^{\alpha}f_{\nu})(x))
H(\lambda (x-y)) \ dy .
$$

We denote the terms on the right-hand side of the last equality by
$I_{1, \alpha}(x)$ and $I_{2, \alpha}(x)$, respectively.
Let
$
K_{\nu, m} = \max \limits_{x \in {\mathbb R}^n,
\vert \beta \vert \le m + 1} \vert (D^{\beta}f_{\nu})(x) \vert $.
Then elementary estimates yield
$$
\max_{x \in {\mathbb R}^n,
\vert \alpha \vert \le m} \vert I_{1, \alpha}(x) \vert \le
\displaystyle \frac  {{\pi}^{\frac n 2} \sqrt n
C_H K_{\nu, m}}
{A \Gamma (\frac n 2 + 1)}
{\lambda}^{-\frac {n}{2n+1}};
$$
$$
\max_{x \in {\mathbb R}^n,
\vert \alpha \vert \le m} \vert I_{2, \alpha}(x) \vert \le
\displaystyle \frac  {2 C_H K_{\nu, m}} {A}
\displaystyle
\int \limits_{\Vert u \Vert > {\lambda}^{\frac 1 {2n +1}}} H(u) \ du .
$$
It follows from these two inequalities that
$$
\displaystyle
\max \limits_{x \in {\mathbb R}^n, \vert \alpha \vert \le m}
\vert
(D^{\alpha}f_{\nu, \lambda})(x) -
(D^{\alpha}f_{\nu})(x)
\vert \to 0
$$
as $\lambda \to +\infty .$
Hence, $q_m (f_{\nu,\lambda} - f_{\nu}) \to 0 $
as $ \lambda \to +\infty $.
Now, since $ m \in {\mathbb N} $ is arbitrary,
it follows that
$ f_{\nu,\lambda} \to f_{\nu} $ in $ {\cal E}(\varphi) $ as
$ \lambda \to +\infty $.

3. We now fix $\lambda>0$ and ${\nu} \in {\mathbb N} $ and approximate
$f_{\nu,\lambda}$ by polynomials in ${\cal E}(\varphi)$.

For $ N \in {\mathbb N}$ let
$$
U_N(x) =
\displaystyle \sum \limits_{k=0}^{N}
\frac
{\sum_{\vert \alpha \vert=k}
(D^{\alpha} H)(0) x^{\alpha}}{k!} \ , \
x \in {\mathbb R}^n.
$$
Since
$$
H(x) = U_N(x) + \displaystyle \sum \limits_{\vert \alpha \vert = N+1}
\displaystyle \frac {(D^{\alpha} H)
(tx) x^{\alpha}}{(N+1)!} \ ,
$$
where $t$ is a point in $(0, 1)$
depending on
$ x \in {\mathbb R}^n $,
it follows with the help of inequality (4) that
\begin{equation}
\vert H(x) - U_N(x) \vert \le
\displaystyle \frac {C_H (N+2)^n {\Vert x \Vert}^{N+1}}{(N+1)!} \ .
\end{equation}

Let $R > 0$ be such that $ supp f_{\nu} \subset  \Pi_R$.
Then we set
$$
V_N(x) =
\displaystyle \frac  {{\lambda}^n} {A}
\int_{\Pi_R}
f_{\nu}(y) U_{N}(\lambda (x-y)) \ d y , \
x \in {\mathbb R}^n.
$$
It is obvious that $V_N$ is a polynomial of degree at most $N$.
We claim that the sequence $(V_N)_{N=1}^{\infty}$ converges to
$ f_{\nu,\lambda} $ in ${\cal E}(\varphi)$ as $N \to \infty $.
Let $m \in {\mathbb N}$ be an arbitrary.
For $\alpha \in \mathbb Z_+^n$ and $x \in {\mathbb R}^n$
$$
(D^{\alpha} f_{\nu,\lambda})(x) - (D^{\alpha} V_{N})(x))
=
\displaystyle \frac {{\lambda}^n} {A}
\int_{\Pi_R}
(D^{\alpha} f_{\nu})(y)
(H({\lambda}(x-y)) - U_N(\lambda (x-y))) \ d y .
$$
Hence with the use of inequality (5) we find from elementary estimates
positive constants $C_1$ and $C_2$ such that for arbitrary $N {\in
{\mathbb  N}}, \vert \alpha \vert \le m, x \in {\mathbb R}^n $
$$
\vert
(D^{\alpha} f_{\nu,\lambda})(x) - (D^{\alpha} V_{N})(x))
\vert
\le \displaystyle \frac {C_1 C_2^N (1+\Vert x \Vert)^{N+1}}
{(N+1)!} \ .
$$
Thus, for each
$N {\in {\mathbb  N}}$
$$
q_m(f_{\nu,\lambda} - V_N) \le
\displaystyle \frac
{C_1 C_2^N}{(N+1)!}
\sup \limits_{x \in {\mathbb R}^n}
\displaystyle \frac {(1+\Vert x \Vert)^{N+1}} {\theta_m(x)} \ .
$$

Let $S^{n-1} = \{x \in {\mathbb R}^n: \Vert x \Vert = 1 \}$. For each
$\sigma \in S^{n-1}$ let $\varphi_{m, \sigma}(t) =
\varphi_{m}(\sigma t), \ t \ge 0$.

Making elementary computations and using Lemma  we get
for some positive constants
$C_3$ and $C_4$ not depending on $N$
$$
q_m(f_{\nu,\lambda} - V_N) \le \displaystyle \frac
{C_3 C_4^N}{(N+1)!}
\displaystyle \frac
{(N+1)^{N+1}}
{\exp(\displaystyle \inf \limits_{\sigma \in S^{n-1}}
(\varphi^*_{m, \sigma}(e^t))^*(N+1))} \ .
$$
Using Stirling formula and the fact that uniformly by $\sigma
\in S^{n-1}$
$$
\displaystyle \lim \limits_{x \to +\infty}
\displaystyle \frac
{(\varphi^*_{m, \sigma}(e^t))^*(N+1)} {N+1}= + \infty
$$
we conclude that
$
q_m(f_{\nu,\lambda} - V_N) \to 0
$
as $N \to \infty $.

The density of polynomials in ${\cal E}(\varphi)$ follows
from  steps 1) -- 3).

2. {\it  On strong dual for ${\cal E}(\varphi)$ in case $n=1$ and
special $\varphi$}. Let $\varphi = \{\varphi_m\}_{m=1}^{\infty}$ be a
family of convex functions $\varphi_m: {\mathbb R} \to {\mathbb R}$
such that:

1. $
\lim \limits_{x \to \infty}
\displaystyle \frac {\varphi_m(x)} {\vert x \vert} = +\infty $
for each $m \in {\mathbb N}$;

2. $\varphi_m(x) - \varphi_{m+1}(x) \ge \ln(1+\vert x \vert) , \
 m \in {\mathbb N}$.

Let
$
\tilde \varphi_m(x)
= \sup \limits_{y \in {\mathbb R}} (x y - \varphi(y)), \
x \in {\mathbb R}.
$

Let
$
Q(\tilde \varphi)
= \bigcup \limits_{m \in {\mathbb N}}
Q(\tilde {\varphi}_m)$,
where
$$
Q(\tilde {\varphi}_m) = \left\{f \in H({\mathbb C}):
N_m(f) = \sup_{z \in  {\mathbb C}}
\displaystyle \frac {\vert f(z) \vert} {(1 + \vert z \vert)^m
\exp(\tilde {\varphi}_m( Im \ z ))} < \infty \right\}.
$$
We endow $Q(\tilde \varphi)$ with the inductive limit topology of the
normed spaces $Q(\tilde {\varphi}_m)$.

For a functional $F \in {\cal E}'(\varphi)$ let
$\hat F(\lambda) = F(e^{-i\lambda x}), \
\lambda \in {\mathbb C}$.

Let ${\cal E}^*(\varphi)$ be the strong dual space for ${\cal
E}(\varphi)$.

Using Theorem 1 and results of [5] (in particular, Theorem 2) we
immediately get

\begin{theorem}
The map
${\cal A}: F \in {\cal E}^*(\varphi) \to \hat F$
establishes topological isomorphism of the spaces
${\cal E}^*(\varphi)$ and $Q(\tilde \varphi)$.

\end{theorem}

3. {\it  On strong dual for a weighted space of sequences of functions}.
On a base of ${\cal E}(\varphi)$ we construct a weighted space of sequences
of infinitely differentiable functions.
Let
$\{(c_k^{(m)})_{k=1}^{\infty}\}_{m=1}^{\infty}$ be a family of sequences
$(c_k^{(m)})_{k=1}^{\infty}$ of positive numbers
$c_k^{(m)}$ such that for each $m \in {\mathbb N}$ \
$
\displaystyle
\sum_{k=1}^{\infty} \displaystyle \frac
{c_k^{(m)}} {c_k^{(m+1)}} = K_m < \infty.
$

Now for each $m \in {\mathbb N}$ we define the normed space
$$
T_m = \left\{f
= (f_1, \ldots, f_k, \ldots ), f_k \in C^m({\mathbb R}):
p_m(f) = \sum_{k=1}^{\infty} c_k^{(m)} q_m(f_k) < \infty \right\}.
$$
Clearly, $T_{m+1} \subset T_m$ for each $m \in {\mathbb N}$.

Let $T = \displaystyle \bigcap \limits_{m=1}^{\infty} T_m$. The linear
space $T$ is equipped with the topology of projective limit of the
spaces $T_m$. Let $T'$ be the space of linear continuous functionals
on $T$, $T^*$ be the strong dual for $T$.

Let $R_m$ be a closure $T$ in $T_m$. $R_m$ is considered with the
topology induced from $T_m$. Then $T$ is a projective limit of the
spaces $R_m$ and $T^*$ is an inductive limit of the spaces $R_m^*$.

Here we describe the space $T^*$ in terms of the Fourier-Laplace
transform of functionals as some weighted space of sequences of entire
functions.

3.1. {\it On the weighted space $T$ of sequences of infinitely
differentiable functions}. It is easy to show that the embeddings
$\i_{m+1, m}: T_{m+1} \to T_m$ are continuous for each $m \in {\mathbb
N}$. Moreover, the embeddings $\i_{m+1, m}: T_{m+1} \to T_m$ are
completely continuous for each $m \in {\mathbb N}$. Really, let $B_r=
\{ f \in T_{m+1}: p_{m+1}(f) \le r \}, r>0,$ be the ball in $T_{m+1}$.
Consider the set ${\cal B}_r$ of all $f = (f_1, \ldots, f_k, \ldots )
\in T_{m+1}$ such that $q_{m+1}(f_k) < (K_m c_k^{(m)})^{-1}, \ k \in
{\mathbb N}$. Obviously, ${\cal B}_r \subset B_r$. As we know for each
$m \in {\mathbb N}$ the space ${\cal E}(\varphi_{m+1})$ is completely
continuously embedded into the space ${\cal E}(\varphi_m)$ [1]. Thus,
the set ${\cal C}_1 = \{u \in {\cal E}(\varphi_{m+1}):q_{m+1}(u) < 1
\}$ is relatively compact in ${\cal E}(\varphi_m)$. So one can extract
a sequence $(u_j)_{j=1}^{\infty}, u_j \in {\cal C}_1,$ convergent in
${\cal E}(\varphi_m)$ to some $u_0 \in {\cal E}(\varphi_m)$ as $j \to
\infty$. Put $f_{j, k}(x) = u_j(x) (K_m c_k^{(m+1)})^{-1}$,
$f_{0,k}(x) = u_0(x) (K_m c_k^{(m+1)})^{-1}, \ j, k \in {\mathbb N}$.
Let $f_j = (f_{j, 1}, \ldots, f_{j, k}, \ldots ), \ j \in {\mathbb
N}$, $f_0 = (f_{0, 1}, \ldots, f_{0, k}, \ldots )$. Then $f_j \in B_r$
for each $j \in {\mathbb N}$ and
$$
p_m(f_j - f_0) =
\displaystyle \sum_{k=1}^{\infty} \displaystyle \frac
{c_k^{(m)} q_m(u_j - u_0)}{K_m c_k^{(m+1)}} \to 0, \ j \to \infty.
$$
Hence, the sequence $(f_j)_{j=1}^{\infty}$ converges in $T_m$ to $f_0$
as $j \to \infty$.  Thus, the space $T$ is an $(M^*)$-space.

3.2. {\it Representation of a linear continuous functional on $T$}.
Let $F \in T'$. Then there exists $m \in {\mathbb N}$, $c>0$ such that
$\vert F(f)\vert \le c p_m(f), \ f \in T$. Define the functional $F_k$
on ${\cal E}(\varphi)$ by the rule $F_k(u) = F(u_k)$, where $u_k$ is
an element of $T$ with $u$ on $k$-th place and zero other components,
$k
\in {\mathbb N}$. We have
\begin{equation}
\vert F_k(u)\vert \le c c_k^{(m)} q_m(u), \
u \in {\cal E}(\varphi).
\end{equation}

Further, for $f = (f_1, \ldots , f_n,  \ldots ) \in T$ \
let
$f_j =(f_1, \ldots , f_j, 0, 0, \ldots ), \ j=1, 2, \ldots$.
Then $f_j \to f$ in $T$
as $j \to \infty$. Consequently, $F(f_j) \to F(f)$
as $j \to \infty$. Thus,
\begin{equation}
F(f) = \
\displaystyle \sum_{j=1}^{\infty} F_j(f_j), \
f \in T.
\end{equation}

3.3. {\it Some notatations and definitions}.
For $k \in {\mathbb N}$, $\lambda \in {\mathbb C}$
\
let $(e^{-i\lambda x})_k$ let be the sequence with
$e^{-i\lambda x}$ on $k$-th place and zero other components,
$x \in {\mathbb R}$.

For a functional $F \in T'$ let $\hat F = ({\hat F_1}, \ldots , {\hat
F_k} , \dots )$, where ${\hat F}_k(\lambda) = F((e^{-i\lambda x})_k),
\
\lambda \in {\mathbb C}, k \in
{\mathbb N}$. Each function ${\hat F_k}$ is an entire function [5].

For each $m \in {\mathbb N}$ let
$$
P_m = \{g = (g_k)_{k=1}^{\infty}, g_k \in H({\mathbb C}):
\Vert g \Vert_m =
$$
$$
=\displaystyle \sup_{k \ge 1} \sup_{\lambda \in {\mathbb C}}
\frac {\vert g_k(\lambda)\vert}
{c_k^{(m)}(1+\vert \lambda \vert)^m \exp (\tilde \varphi_m(Im \
\lambda))} < \infty \}.
$$
Set
$P = \displaystyle \bigcup \limits_{m=1}^{\infty} P_m$.
The linear space $P$ is equipped with the topology of inductive limit of
the spaces $P_m$. It is not difficult to show that $P$ is an
$(LN^*)$-space (for definition of $(LN^*)$-space see [2]).

3.4. {\it Description of $T^*$}. The following result holds.
\begin{theorem} Let functions ${\varphi}_m$ be convex, \ $m \in {\mathbb N}$.
Then the transformation ${\cal F}: S \in T^* \to \hat S$
establishes a topological isomorphism between the spaces $T^*$ and $P$.
\end{theorem}

{\bf Proof}.
Note that if $S \in T^*$ then for some $m \in {\mathbb N}$
$S \in R_m^*$.

Let $F \in R_m^*$, $m \in {\mathbb N}$.
Then
$\vert F(f) \vert \le N_m(F) p_m(f), \ f \in R_m,
$
where $N_m(F)$ is a norm of functional $F$ in $R_m^*$.
For each $k \in {\mathbb N}$
$$
\vert \hat {F_k}(\lambda) \vert \le N_m(F) c_k^{(m)}(1+\vert \lambda \vert)^m
\exp (\varphi^*_m(Im \ \lambda)), \
\lambda \in {\mathbb C}.
$$
Thus, ${\Vert \hat F \Vert}_m \le N_m(F)$, \
$F \in R_m^*$, $m \in {\mathbb N}$.
This means that the linear map ${\cal F}$ acts from $T^*$ to $P$ continuously.

The linear map ${\cal F}$ is injective. For proving it is sufficient
to verify that the system of sequences $(x^{\alpha})_k$ where
$x^{\alpha}$ is on $k$-th place and other components are zero ($\alpha
\in {\mathbb Z}_+, k \in {\mathbb N} $) is complete in $T$. Actually,
let $F \in T'$ be such that $F((x^{\alpha})_k) = 0$ for each $\alpha
\in {\mathbb Z}_+, k \in {\mathbb N} $. Using the representation of
the form (7) we have for a linear continuous functional $F_k$ on
${\cal E}(\varphi)$ that $F_k(x^{\alpha})=0$ for each $\alpha \in
{\mathbb Z}_+$. But then $F_k \equiv 0$, since polynomials are dense
in ${\cal E}(\varphi)$.

The map ${\cal F}$ is surjective.
Let $g = (g_1, \ldots , g_k, \ldots ) \in P$. This means that
$g_k$ are entire functions which satisfy for some $c>0$ the estimate
$$
\vert g_k(\lambda) \vert \le c
c_k^{(m)}
(1+\vert \lambda \vert)^m
\exp(\tilde \varphi_m(Im \ \lambda)), \
k \in {\mathbb N}.
$$
By theorem 2 in [5] (see the proof there, we use the second condition
on functions ${\varphi}_k$) there exist functionals $G_k \in {\cal
E}'(\varphi)$ and a constant $c_1 >0$ such that $G_k(e^{-ix \lambda})
= g_k(\lambda)$ and
$$
\vert G_k(u) \vert \le c_1 c
c_k^{(m)} q_{m+2, 2}(u), \
u \in {\cal E}(\varphi).
$$
We define the functional $G$ on $T$ by the rule:
$$
G(f)=\displaystyle \sum_{k=1}^{\infty} G_k(f_k), \
f = (f_1, \ldots , f_k , \ldots ) \in T.
$$
The series in the right-hand side converges:
$$
\vert G(f) \vert \le c_1 c \displaystyle \sum_{k=1}^{\infty}
c_k^{(m)} q_{m+2}(f_k) = c_1 c p_{m+2}(f), \ f \in T.
$$
Also from this estimate it follows that the linear functional $G$
is continuous. Obviously, $\hat G = g$.

Now the assertion of Theorem 3 follows due the open mapping theorem.

\vspace{1cm}
\begin{center}
REFERENCES
\end{center}

1. V.V. Zharinov, {\it Compact families of LCS and spaces $FS$ and

$DFS$}, Uspekhi matematicheskikh nauk. 34 (1979), pp. 97-131.

2. J. Sebastio e Silva, {\it Su certe classi di spazi localmente
convessi

importanti per le applicazioni}, Rend. Mat. e Appl. 14 (1955).

3. L. Ehrenpreis, {\it Fourier analysis in several complex variables}, New

York: Wiley-Interscience publishers, 1970.

4. B. A. Taylor, {\it Analytically uniform spaces of infinitely differentiable

functions}, Communications on pure and applied mathematics. 24 (1971),

no. 1, pp. 39-51.

5. I. Kh. Musin, {\it Fourier-Laplace transformation of functionals on

a weighted space of infinitely differentiable functions}, Matematicheskii

Sbornik. 191 (2000), no. 10, pp. 57-86.

6. I. Kh. Musin. Fourier-Laplace transformation of functionals on

a weighted space of infinitely differentiable functions on ${\mathbb R}^n$.

Matematicheskii Sbornik, 195 (2004),  no. 10, pp. 83-108.

7. P. Koosis. {\it Introduction to $H_p$ spaces}, Cambridge Univ.
Press.

Cambridge. 1980

\vspace {1cm}
Institute of mathematics with computer center of Ufa scientific center of
Russian Academy of Sciences,
Chernyshevskii st., 112, Ufa, Russia

e-mail: musin@imat.rb.ru
\end{document}